\documentclass[reqno, 11pt]{amsart}

\usepackage[letterpaper,hmargin=1in,vmargin=1.2in]{geometry}
\usepackage{amsmath, amssymb, amsthm, verbatim, url}
\usepackage{graphicx}
\usepackage{enumerate, pinlabel}
\usepackage{enumitem}
\usepackage{amsmath,amssymb,amsthm,mathrsfs,graphicx,url}
\usepackage[usenames,dvipsnames]{color}\usepackage[colorlinks=true,linkcolor=Black,citecolor=Black]{hyperref}

\newtheorem*{theorem}{Theorem}
\title[Manifolds having a finite and positive number of embedded eigenvalues]{Manifolds with cylindrical ends having a finite and  positive number of embedded eigenvalues}
\author{T. J. Christiansen and K. Datchev}
\address{Department of Mathematics, University of Missouri, Columbia, MO 65211 USA}
\email{christiansent@missouri.edu}
\address{Department of Mathematics, Purdue University, West Lafayette, IN 47907 USA}
\email{kdatchev@purdue.edu}

\begin{document}

\begin{abstract}
We construct a surface with a cylindrical end which has a finite number of Laplace eigenvalues embedded in its continuous spectrum. The surface is obtained by attaching a cylindrical end to a hyperbolic torus with a hole.  To our knowledge, this is the first example of a manifold with a cylindrical end whose number of 
eigenvalues is known to be finite and nonzero. The construction can be varied to give examples with arbitrary genus and with an arbitrarily large finite number of 
eigenvalues. The constructed surfaces also have resonance-free regions near the continuous spectrum and long-time asymptotic expansions of solutions to the wave equation.
\end{abstract}

\maketitle

Let $(X,g)$ be a manifold with cylindrical ends and let $-\Delta \ge 0$ be its Laplacian. A basic question is: How many eigenvalues does $-\Delta$ have?
This question is subtle because the essential spectrum of $-\Delta$ is $[0,\infty)$, and hence any eigenvalues must be embedded in it rather than isolated. 
In some examples, including the simplest example of the cylinder $\mathbb R \times (\mathbb R / \mathbb Z)$,  the number of eigenvalues is zero. In this note we prove that the number of eigenvalues can be nonzero but finite. 

\begin{theorem}
There exists a surface with a cylindrical end whose Laplacian has a nonzero but finite number of eigenvalues.
\end{theorem}

\begin{figure}[h]
\labellist
\small
\pinlabel $\ell$ at 98 36
\pinlabel $\ell$ at 281 36
\endlabellist
\includegraphics[height=3cm]{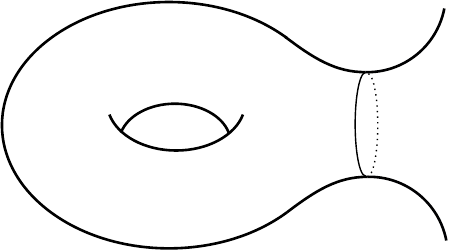}
\hspace{2cm}
\includegraphics[height=3cm]{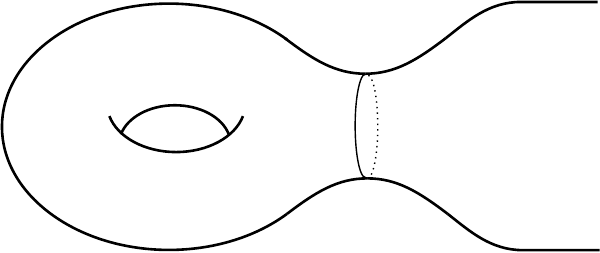}
\caption{The hyperbolic surface $(X,g_H)$, and the surface with a cylindrical end $(X,g)$.
}\label{f:surfaces}
\end{figure}

In fact, we construct families of Riemannian surfaces $(X,g)$ consisting of a compact hyperbolic piece with an infinite cylindrical end attached. More specifically, we start with a symmetric funneled torus of constant negative curvature, $(X,g_H)$, and then modify the metric in the funnel end to obtain a cylindrical end (a Riemannian product of a half-line with a circle) as in Figure~\ref{f:surfaces}.

In the proof below, we show that if the length $\ell$ of the geodesic at the boundary of the funnel is  sufficiently small, then $-\Delta$ has an eigenvalue. For this, as in work of earlier authors including \cite{elv, dp}, we consider  $-\Delta|_O$, a restriction of the Laplacian $-\Delta$ to   a domain of functions with a suitable symmetry, and we show that $-\Delta|_O$ has an isolated eigenvalue. 

The more difficult step is showing that $-\Delta$ does not have infinitely many eigenvalues. This is done in \cite[Theorem 1.1]{cdaens} for a family of manifolds with cylindrical ends which includes $(X,g)$, and which also includes the variants of $(X,g)$ described in the extensions after the proof below. This result is an application of two recent breakthroughs in the scattering theory of hyperbolic manifolds: Vasy's propagation of semiclassical singularies \cite{vasy}, and the essential spectral gap of Dyatlov and Zahl \cite{dyza} and Bourgain and Dyatlov \cite{bd}.

Our results also imply the existence of resonance-free regions and asymptotic expansions for solutions to the wave equation $(\partial_t^2 - \Delta)u=0$. See the discussion following the proof below for more.

Existence proofs for  eigenvalues for manifolds with cylindrical ends (including the case of waveguides, or more generally manifolds with boundary, where we may have isolated eigenvalues below the essential spectrum as well as embedded ones within it)  have been given by many authors under various assumptions. In some cases the number of eigenvalues is shown to be infinite \cite{w,cz,parn, ik}, and in many others it is left open whether the number of embedded eigenvalues is finite or infinite: see \cite{es, abgm, goja, elv, bgrs, dp} for some such results, 
\cite{lcm,kk,ek} for further results and surveys of the literature, and \cite{cdn,sm} for some recent results and more references. In general, the eigenvalue counting function is 
polynomially bounded \cite{d,cz,parn}. On the negative side, nonexistence results for eigenvalues have been obtained in many of the aforementioned papers, as well as in \cite{MoWe,bdk,cdstar}, and one for a resonance at the bottom of the spectrum in \cite{gj}. To our knowledge, the examples in this note are the first manifolds with cylindrical ends whose total number of eigenvalues 
is known to be finite but not zero.

We discuss briefly other geometric situations in which finite numbers of embedded eigenvalues  arise. Particularly important, and closely related to manifolds with cylindrical ends \cite{d,gui},  are manifolds with cusps: 
see \cite{cdv}, the surveys \cite{sarnak, muller}, and  \cite{bon, hj} for  recent results and pointers to the substantial literature on this topic. Another direction of study is manifolds whose curvature converges to a constant at infinity. One asks: For what rates of convergence do embedded eigenvalues occur? See \cite{d10,jili} for some recent results and references, and \cite{ku} for an example with precisely one embedded eigenvalue.

\begin{proof}[Proof of Theorem]
We construct the desired Riemannian surface $(X,g)$ by modifying the metric of a certain hyperbolic surface $(X,g_H)$, which is the quotient of the hyperbolic plane $\mathbb H$ by a Schottky group. We use  the unit disk model $\mathbb H = \{z \in \mathbb C \mid |z| <1\}$. For each $j \in \{1, \, 2, \, 3, \, 4\}$, let 
$D_j$ be the open disk in $\mathbb C$, orthogonal to the unit disk, which is inscribed in the sector $\{z \in \mathbb C \mid |\arg z - (2j-1)\pi/4| < \alpha\}$, 
where $\alpha < \pi/4$ will be chosen later (close to $\pi/4$). See Figure~\ref{f:disk}.  

\begin{figure}[h]
\labellist
\small
\pinlabel $D_1$ at 270 270
\pinlabel $D_2$ at 100 270
\pinlabel $D_3$ at 100 100
\pinlabel $D_4$ at 270 100
\pinlabel $\mathcal F$ at 152 152
\pinlabel $B$ at 705 177
\pinlabel $D_1'$ at 705 315
\pinlabel $D_2'$ at 560 177
\pinlabel $D_3'$ at 705 40
\pinlabel $D_4'$ at 850 177
\scriptsize
\pinlabel $2\alpha$ at 203 197
\tiny
\pinlabel $a$ at 820 215
\pinlabel $a$ at 666 65
\pinlabel $b$ at 740 65
\pinlabel $b$ at 588 216
\pinlabel $c$ at 740 293
\pinlabel $c$ at 588 145
\pinlabel $d$ at 666 295
\pinlabel $d$ at 820 145
\endlabellist
\includegraphics[width=5cm]{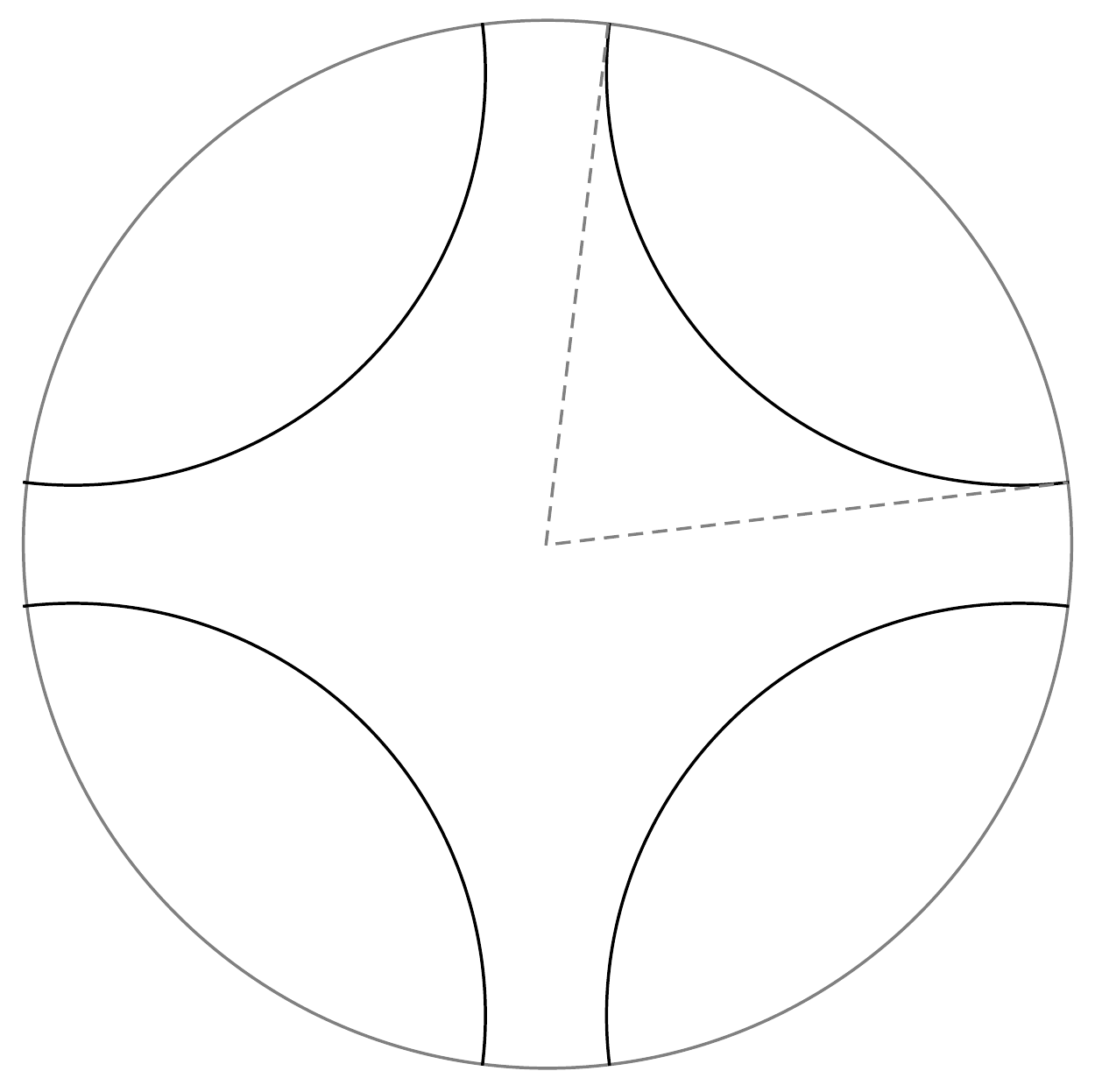}
\hspace{2cm}
\includegraphics[width=5cm]{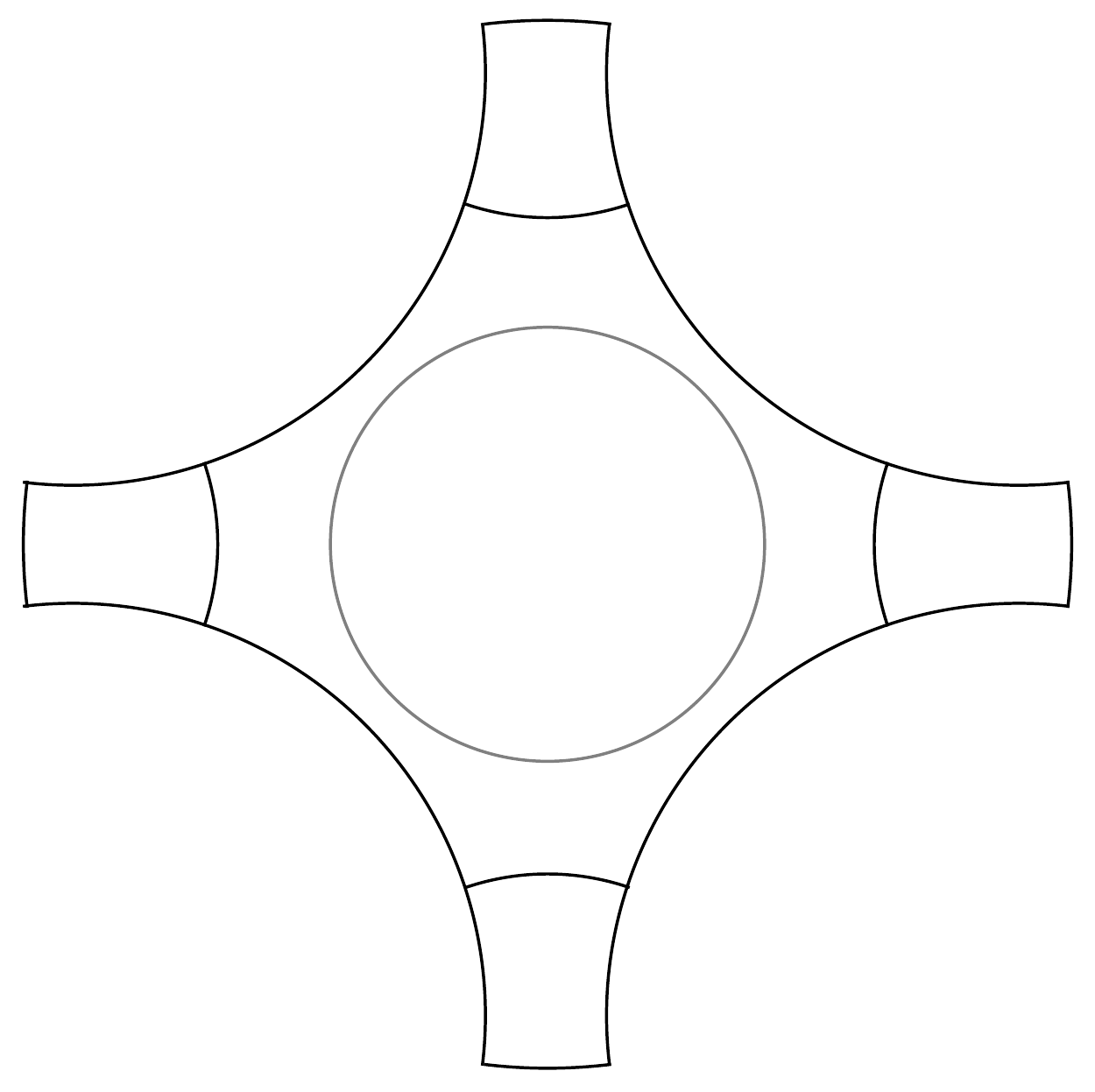}
\caption{
On the left, the unit disk with the  parts of the disks $D_1$, $D_2$, $D_3$, and $D_4$ used to define the fundamental domain $\mathcal F$. On the right, the fundamental domain $\mathcal F$, with points having matching letters identified.
}\label{f:disk}
\end{figure}
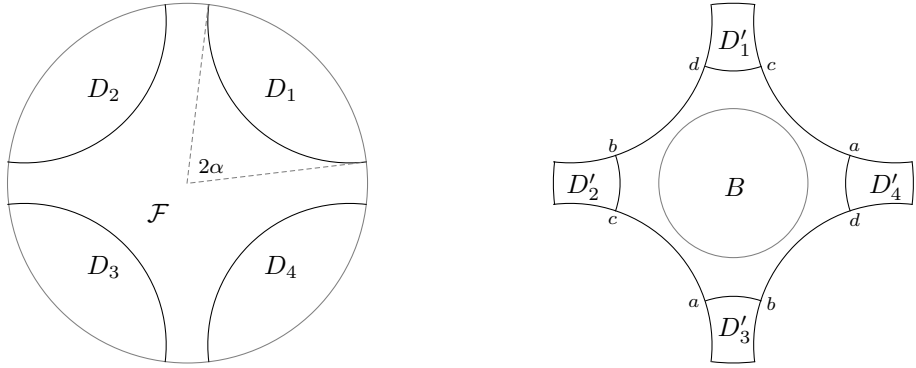

Let $S_1$ be the orientation-preserving isometry of $\mathbb H$ which fixes $e^{\pi i /4}$ and $-e^{ \pi i/4}$ and which maps $D_1$ to the exterior of $D_3$: see  \cite[Section~2.1]{bo} and \cite[Project~4.2]{msw}.
Let 
\begin{equation}\label{e:s1s2bar}
S_2(z) = \overline{S_1^{-1}(\bar z)}.
\end{equation}
Then $S_2$ maps  $D_2$ to the exterior of $D_4$.

Let  $\Gamma$ be the group generated by $S_1$ and $S_2$, and let $(X,g_H)$ be the quotient of $\mathbb H$ by $\Gamma$. We recall some standard facts about $(X,g_H)$ from \cite[Sections~2.4 and~15.1]{bo}. The surface $(X,g_H)$ is obtained from the fundamental domain
\[
\mathcal F:=\mathbb H -\bigcup_{j=1}^{4}D_j
\]
by identifying $\partial D_1$ with $\partial D_3$ using $S_1$ and $\partial D_2$ with $\partial D_4$ using $S_2$, so that points having matching letters in Figure \ref{f:disk} are identified. Let $N$ be the convex core of $X$; this is the part of $X$ whose lift to $\mathcal F$ is given by $\mathcal F -  \bigcup_{j=1}^{4}D_j'$, where each $D_j'$ is the open disk in $\mathbb C$, orthogonal to the unit disk,
to $D_j$, and to $D_{j+1}$ (with the convention that $D_5=D_1$),  whose center has argument $j\pi/2$. 
Then $N$ has genus $1$
and  $\partial N$ is a geodesic of length  $\ell = 4 \rho(D_1,D_2)$,  where $\rho$ is the hyperbolic distance function.  Furthermore, $X \setminus N$ is a hyperbolic funnel, i.e.  there are coordinates $(r,y)\in(0,\infty) \times (\mathbb R / \ell \mathbb Z)$ such that
\begin{equation}\label{e:end}
 X\setminus N = (0,\infty)_r\times (\mathbb R / \ell \mathbb Z)_y, \qquad g_H|_{X \setminus N} = dr^2 + \cosh^2 \!r\, dy^2.
\end{equation}
By \eqref{e:s1s2bar}, the involution $z \mapsto \bar z$ on $\mathcal F$ descends to an involution  $J$ on $X$, and $J|_{ X\setminus N }$ is a reflection. For later convenience, we normalize the coordinate $y$ in \eqref{e:end} so that $J(r,y) = (r,-y)$. For example, we may take $y=0$ when $\arg z = 0$, $y = \ell/8$ at $a$, $y=3\ell/8$ at $b$, $y= \pm \ell/2$ when $\arg z = \pm \pi$, $y= -3\ell/8$ at $c$, $y = -\ell/8$ at $d$.

Let $B= \{z \in \mathbb C \mid |z| <\sqrt 2 - 1\}$, and fix $\varphi \in C_c^\infty(B)$ such that $\varphi \not\equiv 0$ and $\varphi(\bar z) = - \varphi( z)$ for all $z$. Since $B \subset \mathcal F$, regardless of $\alpha$, we may identify $\varphi$ with a function on $X$, and $\varphi \circ J= - \varphi$. Moreover, this $\varphi$ is supported in $N$, regardless of $\alpha$.

Now choose $\alpha<\pi/4$ (sufficiently close to $\pi/4$) such that $\ell = \ell(\alpha)$ obeys
\begin{equation}\label{e:elle}
\frac{\|d\varphi\|^2_{L^2(X,g_H)}}{\|\varphi\|^2_{L^2(X,g_H)}}<\frac 1 {\ell^2}.
\end{equation}
More precisely, $\tan(\alpha) \sinh  \rho(0,D_1)  =1$ by the angle of parallelism formula \cite[Theorem~7.9.1]{be}, and $\cosh(\ell/4)  = \sinh^2 \! \rho(0,D_1)$ by a property of right-angled pentagons  \cite[(7.18.2)]{be}. See Figure~\ref{f:hyptrig}. Hence, $\tan^2(\alpha) \cosh(\ell/4)=1$, and  \eqref{e:elle} is equivalent to $\tan^2(\alpha)\cosh (\|\varphi\|/4\|d\varphi\|) >1$.

\newcommand*{\vcenteredhbox}[1]{\begin{tabular}{@{}c@{}}#1\end{tabular}}

\begin{figure}[h]
\labellist
\small
\pinlabel $\ell/4$ at 735 70
\footnotesize
\pinlabel $\rho(0,D_1)$ at 28 85
\pinlabel $\rho(0,D_1)$ at 525 135
\scriptsize
\pinlabel $\alpha$ at 40 23
\endlabellist
\vcenteredhbox{\includegraphics[width=5cm]{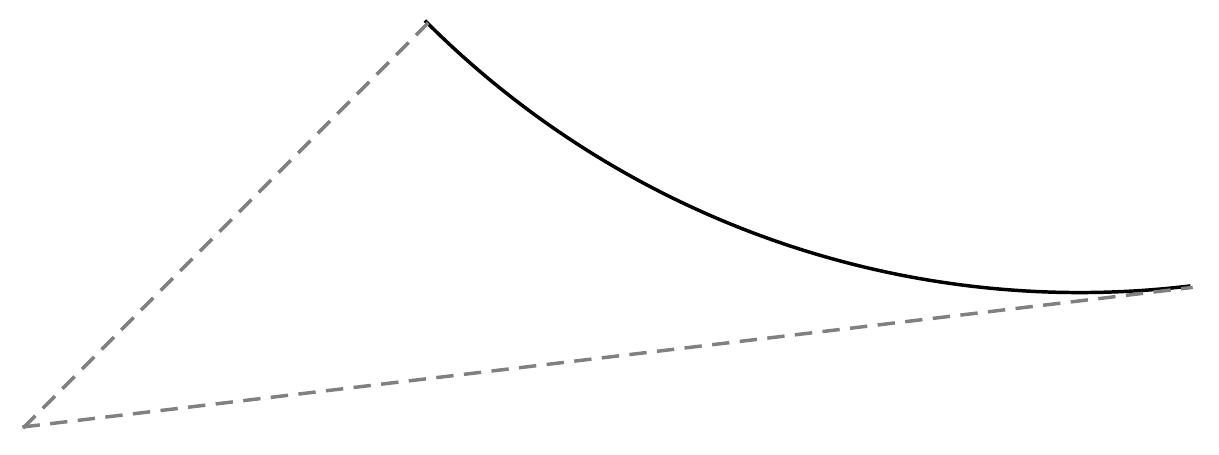}}
\hspace{2cm}
\vcenteredhbox{\includegraphics[width=3cm]{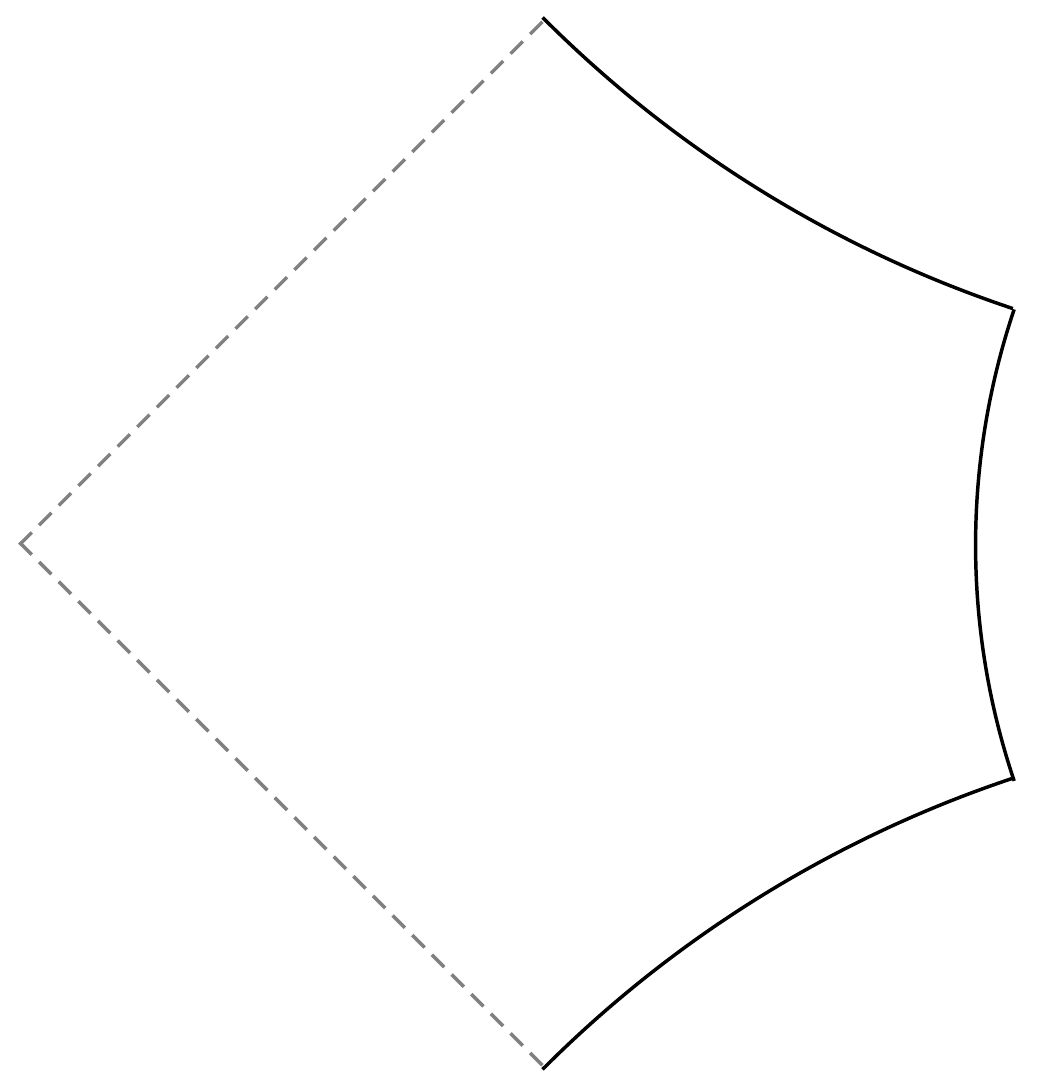}}
\caption{The triangle and pentagon used to obtain \eqref{e:elle}.
}\label{f:hyptrig}
\end{figure}

Next, modify the metric in the funnel end so as to change it into a cylindrical end in the following way.  Fix a metric $g$ on $X$ having the form
\[
g|_N = g_H|_N, \qquad g|_{X \setminus N} = dr^2 + F(r)dy^2,
\] 
where $F \in C^\infty((0,\infty))$ obeys $F(r) = \cosh^2\!r$ near $r=0$, $F'(r)>0$ when $ r \in (0,R)$ for some $R>0$, and $F(r) = (2\pi)^2$ for $r \ge R$. Note that $g$, like $g_H$, is invariant under $J$.

Let $-\Delta|_O$ be the Laplacian on $(X,g)$, restricted to functions $u$ obeying $u \circ J = -u$. Observe that we have, in the part of $X\setminus N$ where $r > R$, the Fourier series expansion
\[
-\Delta u(r,y) =\left(-\partial_r^2 - \partial_y^2/(2\pi)^2 \right) \sum_{n =1}^\infty  u_n(r) \sin\left(2\pi ny/\ell\right) 
\ge \left(-\partial_r^2 + 1/\ell^2\right)u(r,y),
\]
for all  $u$ obeying $u \circ J = -u$, with equality if $u_n=0$ for $n \ge 2$. This implies that the essential spectrum of $-\Delta|_O$ is $[1/\ell^2,\infty)$. But, by  \eqref{e:elle}, $\varphi$ is not in the range of the spectral projection of $-\Delta|_O$ onto $[1/\ell^{2},\infty)$.  Hence  $-\Delta|_O$ has at least one eigenvalue, and in particular $-\Delta$ has at least one eigenvalue.

By \cite[Theorem~1.1]{cdaens}, $-\Delta$ does not have infinitely many eigenvalues.
\end{proof}

\subsection*{Extensions} 

\begin{enumerate}[label=(\roman*)]
\item To construct examples with at least $m$  eigenvalues, replace \eqref{e:elle} with
\[
\max_{\varphi \in V} \frac{\|d\varphi\|^2_{L^2(X,g_H)}}{\|\varphi\|^2_{L^2(X,g_H)}}<\frac 1 {\ell^2},
\]
where $V$ is an $m$-dimensional space of functions in $C_c^\infty(B)$ satisfying $\varphi(\bar z) = - \varphi(z)$.

\item The symmetry with respect to the vertical axis (but not that with respect to the horizontal axis) can be broken: after choosing $\alpha$ such that \eqref{e:elle} holds, replace $D_2$ and $D_3$ by disks, orthogonal to the unit disk, inscribed in larger sectors  $\{z \in \mathbb C \mid |\arg z - (2j-1)\pi/4| < \alpha'\}$ for some $\alpha < \alpha' <\pi/4$.

\item To construct examples having genus $g$ for any $g \ge 1$, start with disks $D_1, \, \dots \, D_{4g}$, orthogonal to the unit disk, inscribed in sectors $\{z \in \mathbb C \mid |\arg z - (2j-1)\pi/4g| < \alpha\}$, with $\alpha<\pi/4g$, and take $S_j$ mapping $D_j$ to the exterior of $D_{j + 2g}$. Then the convex core of the quotient has genus $g$, and the boundary is a geodesic of length $\ell = 4g \rho(D_1,D_2)$, and we can choose $\alpha< \pi/4g$, close to $\pi/4g$, such that \eqref{e:elle} holds and proceed as before.

\item In all of these examples, we can also rule out resonances close to the continuous spectrum. We discuss this briefly here, referring the reader to \cite[Theorem 5.6]{cdaens} for details. By \cite[Theorem~1.1]{cdaens}, the cut-off resolvent of $-\Delta$ is bounded at high energies, and by a resolvent identity of Vodev \cite{v}, this implies that there is a polynomial neighborhood of the continuous spectrum in a suitable Riemann surface (see \cite{gui, mel})  to which the cut-off resolvent continues holomorphically.

\item By \cite[Theorem 1.2]{cdwaveexp}, we also have an asymptotic expansion for solutions to the wave equation on $(X,g)$. Let $u(t,x)$ solve
\[
(\partial_t^2 - \Delta ) u(t,x) = 0, \quad
u(0,x) = f_1(x), \quad \partial_t u(0,x) = f_2(x),
\]
where $f_1$ and $f_2$ are in $C_c^\infty(X)$. Then, for any $k_0 \ge 1$, uniformly as $t \to \infty$ and as $x$ varies in a compact subset of $X$,

\[
\begin{split}
 u(t,x)  = &\int_X f_2 + \sum_{m=1}^M  \left( \cos(E_m t) \int_X f_1 \eta_m+ 
 \frac{\sin( E_m t)}{E_m } \int_X f_2\eta_m \right)\eta_m(x) + \\& \sum_{k=0}^{k_0-1}t^{-1/2-k} 
\sum_{n=1}^\infty \Big(\cos(nt/\ell)a_{n,k}(x)+ \sin(nt/\ell)b_{n,k}(x)\Big) + 
\mathcal{O}(t^{-k_0}),
\end{split}
\]
where $\{\eta_1, \dots,\eta_M\}$ is a maximal, real-valued, orthonormal set of eigenfunctions of $-\Delta$,  $-\Delta \eta_m = E_m^2 \eta_m$, and
$a_{n,k},\,b_{n,k}\in C^\infty(X)$. 
By the Theorem, the sum over $m$ is not empty.

\end{enumerate}

\subsection*{Acknowledgments}  The authors gratefully acknowledge the 
partial support of the Simons Foundation (TC, collaboration grant for
mathematicians) and the National Science Foundation (KD,  Grant DMS-1708511). TC is also grateful for partial support from a University of Missouri research leave.
Thanks also to David Borthwick for helpful conversations about  hyperbolic quotients.

\end{document}